\documentclass[11pt]{amsproc} 
\usepackage{graphicx}
\usepackage{amsfonts}
\usepackage{amscd}
\usepackage{amssymb}
\usepackage{alltt}
\usepackage{url}

\usepackage{xcolor}
\usepackage{tikz}
\usetikzlibrary{chains,shapes,arrows,trees,matrix,positioning,decorations,backgrounds,fit}

\newtheorem{thm}{Theorem}

\parskip=\baselineskip

\def\op#1{{\operatorname{#1}}}
\newcommand{\ring}[1]{\mathbb{#1}}
\def\tikzfig#1#2#3{%
\begin{figure}[htb]%
  \centering
\begin{tikzpicture}#3
\end{tikzpicture}
  \caption{#2}
  \label{fig:#1}%
\end{figure}%
}

\def\gg{\mathfrak{g}}
\def\cc{\mathfrak{c}}
\def\AA{{\mathcal A}}
\def\MM{{\mathcal M}}
\def\OO{{\mathcal O}}
\def\UU{{\mathcal U}}

\title{The Work of Ng\^o Bao Ch\^au}
\author{Thomas C. Hales} 
\email{hales@pitt.edu}

\begin{document}

\maketitle

In August 2010, Ng\^o Bao Ch\^au was awarded a Fields Medal for his
deep work relating the Hitchin fibration to the Arthur-Selberg trace
formula, and in particular for his proof of the Fundamental Lemma for
Lie algebras~\cite{NBC:2006},~\cite{NBC:2010}.

\section{the Trace Formula}

A function $h:G\to\ring{C}$ on a finite group $G$ is a {\it class} function
if $h(g^{-1} x g) = h(x)$ for all $x,g\in G$.  A class
function is constant on each conjugacy class.  A basis of the
vector space of class functions is the set of characteristic functions
of conjugacy classes.

A {\it representation} of $G$ is a homomorphism $\pi:G\to GL(V)$ from $G$
to a group of invertible linear transformations on a complex vector
space $V$.  It follows from the matrix identity 
\[
\op{trace}( B^{-1} A B) = \op{trace} (A)
\]
that the function $g\mapsto \op{trace}(\pi(g))$
is a class function.  This function is called an
{\it irreducible character} if $V$ has no proper $G$-stable subspace.
A basic theorem in finite group theory asserts that the set of
irreducible characters forms a second basis of the vector space of
class functions on $G$. 

A {\it trace formula} is an equation that gives the expansion of a
class function $h$ on one side of the equation in the basis of
characteristic functions of conjugacy classes $C$ and on the other
side in the basis of irreducible characters
\[
\sum_C b_C\ \op{char}_C     =  h = \sum_\pi a_\pi\ \op{trace}\, \pi,
\]
for some complex coefficients $b_C$ and $a_\pi$ depending on $h$.
The side of the equation with conjugacy classes is called the {\it
  geometric side} of the trace formula and the side with irreducible
characters is called the {\it spectral side}.

When $G$ is no longer assumed to be finite, some analysis is
required.  We allow $G$ to be a Lie group or more generally a locally
compact topological group.  The vector space $V$ may be 
infinite-dimensional so that a  trace of a linear
transformation of $V$ need not converge.  To improve convergence,
the irreducible character is no longer viewed as a function, but rather as a
distribution
\[
f \mapsto \op{trace} \int_{G} \pi(g) f(g) \, dg,
\]
where $f$ runs over smooth compactly supported test functions on
the group, and $dg$ is a $G$-invariant measure.  Similarly, the
characteristic function of conjugacy class is replaced with a
distribution that integrates a test function $f$ over the
conjugacy class $C$ with respect to an invariant measure:
\begin{equation}\label{eqn:oi}
f \mapsto \int_C f(g^{-1} x g) \, dg.
\end{equation}
The integral (\ref{eqn:oi}) is called an {\it orbital integral}.  A
trace formula in this setting becomes an identity that expresses a
class distribution (called an {\it invariant distribution}) on the
geometric side of the equation as a sum of orbital integrals and on
the spectral side of the equation as a sum of distribution characters.

The celebrated Selberg trace formula is an identity of this general
form for the invariant distribution associated with the representation
of $SL_2(\ring{R})$ on $L^2(SL_2(\ring{R})/\Gamma)$, for a discrete
subgroup $\Gamma$.  Arthur generalized the Selberg trace formula to
reductive groups of higher rank.

\section{History}

The Fundamental Lemma (FL) is a collection of identities of orbital
integrals that arise in connection with a trace formula.  It takes
several pages to write all of the definitions that are needed for a precise
statement of the lemma~\cite{Hales:FL-statement}.  Fortunately, the
significance of the lemma and the main ideas of the proof can be
appreciated without the precise statement.

Langlands conjectured these identities in lectures on the trace
formula in Paris in 1980 and later put them in more precise form with
Shelstad \cite{Langlands:debuts},~\cite{LS:1987}.  Over time,
supplementary conjectures were formulated, including a {\it twisted}
conjecture by Kottwitz and Shelstad, and a {\it weighted} conjecture
by Arthur~\cite{KS:1999},~\cite{Arthur:2002}.  Identities of orbital
integrals on the group can be reduced to slightly easier identities on
the Lie algebra~\cite{LS:1990}.  Papers by Waldspurger rework the
conjectures into the form eventually used by Ng\^o in his
solution~\cite{Wald:2008}, \cite{Wald:1991}.  Over the years,
Chaudouard, Goresky, Kottwitz, Laumon, MacPherson, and Waldspurger
among others have made fundamental contributions that led up to the
proof of the FL or extended the results afterwards~\cite{LN:08},
\cite{GKM:2004}, \cite{GKM:2006}, \cite{CL:2009:I}, \cite{CL:2009:II},
\cite{CHL:2010}.  It is hard to do justice to all those who have
contributed to a problem that has been intensively studied for
decades, while giving special emphasis to the spectacular
breakthroughs by Ng\^o.

With the exception of the FL for the special linear group $SL(n)$,
which can be solved with representation theory; starting in the early
1980s, all plausible lines of attack on the general problem have been
geometric.  Indeed, a geometric approach is suggested by direct 
computations of these integrals in special cases, which
give their values as the number of points on hyperelliptic curves
over finite fields~\cite{KL:1988}, \cite{Hales:hyperelliptic}.

To motivate the FL, we must recall the bare outlines of the ambitious
program launched by Langlands in the late 1960s to use representation
theory to understand vast tracts of number theory.  Let $F$ be a
finite field extension of the field of rational numbers $\ring{Q}$.
The ring of adeles $\ring{A}$ of $F$ is a locally compact topological
ring that contains $F$ and has the property that $F$ embeds discretely
in $\ring{A}$ with a compact quotient $F\backslash\ring{A}$.  The ring
of adeles is a convenient starting point for the analytic treatment of
the number field $F$.  If $G$ is a reductive group defined over $F$
with center $Z$, then $G(F)$ is a discrete subgroup of $G(\ring{A})$
and the quotient $G(F)Z(\ring{A})\backslash G(\ring{A})$ has finite
volume.  A representation $\pi$ of $G(\ring{A})$ that appears in the
spectral decomposition of
\[
L^2(G(F)Z(\ring{A})\backslash G(\ring{A}))
\]
is said to be an {\it automorphic representation}.  The automorphic
representations (by descending to the quotient by $G(F)$) are those
that encode the number-theoretic properties of the field $F$.
The theory of automorphic representations just for the two linear
groups $G=GL(2)$ and $GL(1)$ already 
encompasses the classical theory
of modular forms and global class field theory.

There is a complex-valued
function $L(\pi,s)$, $s\in \ring{C}$, called an automorphic $L$-function,
attached to each automorphic representation $\pi$. (The $L$-function
also depends on a representation of a dual group, but we skip these details.)
Langlands's philosophy can be summarized as two objectives:
\begin{enumerate}
\item Show that many $L$-functions that routinely arise in
number theory are automorphic.
\item Show that automorphic $L$-functions have wonderful analytic properties. 
\end{enumerate}

There are two famous examples of this philosophy.  
In Riemann's  paper on the zeta function
\[
\zeta(s) = \sum_{n=1}^\infty \frac{1}{n^s},
\]
he proved that it has a functional equation and meromorphic
continuation by relating it to a $\theta$-series (an automorphic
entity) and then using the analytic properties of the $\theta$-series.
Wiles proved Fermat's Last Theorem by showing that the $L$-function
$L(E,s)$ of every semi-stable elliptic curve over $\ring{Q}$ is
automorphic.  From automorphicity follows the analytic continuation
and functional equation of $L(E,s)$.

The Arthur-Selberg trace formula has emerged as a general tool to
reach the first objective (1) of Langlands's philosophy.  To relate
one $L$-function to another, two trace formulas are used in tandem
(Figure~\ref{fig:identities}).  An automorphic $L$-function can be
encoded on the spectral side of the Arthur-Selberg trace formula.  A
second $L$-function is encoded on the spectral side of a second trace
formula of a possibly different kind, such as a topological trace
formula.  By equating the geometric sides of the two trace formulas,
identities of orbital integrals yield identities of $L$-functions.
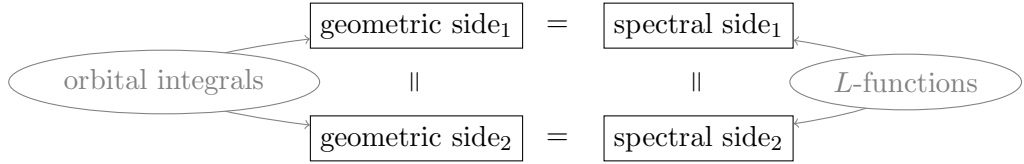
\tikzfig{identities} {A pair of trace formulas can transform
  identities of orbital integrals into identities of $L$-functions.}
{
\begin{scope}[scale=0.25]
\path  (0,0) node(nw)[rectangle,draw,color=black] {geometric side$_1$} -- 
   node[color=black]{$=$} (15,0) node(ne)[rectangle,draw,color=black] {spectral side$_1$};
\path (0,-6) node(sw)[rectangle,draw,color=black] {geometric side$_2$} -- 
  node[color=black]{$=$} (15,-6) node(se)[rectangle,draw,color=black] {spectral side$_2$};
\draw[<->,color=gray] (sw) .. controls (-16,-4) and (-16,-2)  
  .. node[ellipse,draw,fill=white!20] {orbital integrals} (nw);
\draw[<->,color=gray] (se) .. controls (28,-4) and (28,-2) 
   .. node[ellipse,draw,fill=white!20] {$L$-functions} (ne);
\node[rotate=90] at (15,-3) {$=$};
\node[rotate=90] at (0,-3) {$=$};
\end{scope}
}

The value of the FL lies in its utility.  The FL can be characterized
as the minimal set of identities that must be proved in order to put
the trace formula in a useable form for applications to number theory, such as those
mentioned at the end of this report.

\section{the Hitchin Fibration}

Ng\^o's proof of the FL is based on the Hitchin fibration~\cite{Hitchin:87}.

Every endomorphism $A$
of a finite dimensional vector space $V$ has a characteristic polynomial
\begin{equation}
\det(t - A) = t^n + a_1 t^{n-1} + \cdots +a_n.
\end{equation}
Its coefficients $a_i$ are symmetric polynomials of the eigenvalues of
$A$.  This determines a {\it characteristic map} $\chi:\op{end}(V) \to
\cc$, from the Lie algebra of endomorphisms of $V$ to the vector space
$\cc$ of coefficients $(a_1,\ldots,a_n)$.  This construction
generalizes to a characteristic map $\chi:\gg\to \cc$ for every
reductive Lie algebra $\gg$, by evaluating a set of symmetric
polynomials on $\gg$.

Fix once and for all a smooth projective
curve $X$ of genus $g$ over a finite field $k$.  

In its simplest form, a {\it Higgs pair} $(E,\phi)$
is what we obtain when we allow an element $Z$ of the Lie
algebra $\op{end}(V)$ to vary continuously along the curve $X$.  
As we vary along the curve, the vector space $V$ sweeps out a vector
bundle $E$ on $X$, and the element $Z\in\op{end}(V)$ sweeps out a
section $\phi$ of the bundle $\op{end}(E)$ or of the bundle
$\op{end}(E)\otimes {\OO}_X(D)$ when the section acquires finitely
many poles prescribed by a divisor $D$ of $X$.  Extending this
construction to a general reductive Lie group $G$ with Lie algebra
$\gg$, a Higgs pair $(E,\phi)$ consists of a principal $G$-bundle $E$
and a section $\phi$ of the bundle $\op{ad}(E)\otimes {\OO}_X(D)$
associated with $E$ and the adjoint representation of $G$ on $\gg$.
For each $X$, $G$, and $D$, there is a moduli space ${\MM}$ (or more
correctly, moduli stack) of all Higgs pairs $(E,\phi)$.

The {\it Hitchin fibration} is the morphism obtained when we vary the
characteristic map $\chi:\gg\to \cc$  along a
curve $X$.  For each Higgs pair $(E,\phi)$, we evaluate the
characteristic map $p\mapsto \chi(\phi_p)$ of the endomorphism $\phi$
at each point $p\in X$. This function belongs to the set $\AA$ of a
global sections of the bundle $\cc\otimes {\OO}_X(D)$ over $X$.  The
Hitchin fibration is this morphism ${\MM}\to {\AA}$.

{\it Abelian varieties} occur naturally in the Hitchin
fibration.  To illustrate, we return to the Lie
algebra $\gg=\op{end}(V)$.  For each section $a=(a_1,\ldots,a_n)\in
\AA$,  the characteristic
polynomial 
\begin{equation}\label{eqn:spectral}
t^n + a_1(p) t^{n-1} + \cdots+ a_n(p)=0,
\end{equation} 
defines an $n$-fold cover $Y_a$ of $X$ (called the {\it spectral
  curve}).  By construction, each point of the spectral curve is a
root of the characteristic polynomial at some $p\in X$.  We consider
the simple setting when $Y_a$ is smooth and the discriminant of the
characteristic polynomial is sufficiently generic.  A Higgs pair $(E,\phi)$ over
the section $a$ determines a line (a one-dimensional eigenspace of
$\phi$ with eigenvalue that root) at each point of the spectral curve,
and hence a line bundle on $Y_a$.  This establishes a map from Higgs
pairs over $a$ to $\op{Pic}(Y_a)$, the group of line bundles on the
spectral curve $Y_a$.  Conversely, just as linear maps can be
constructed from eigenvalues and eigenspaces, Higgs pairs can be
constructed from line bundles on the spectral curve $Y_a$.  The
connected component $\op{Pic}^0(Y_a)$ is an abelian variety.
Even outside this simple setting,  the group of symmetries of
the Hitchin fiber over $a\in\AA$ has an abelian
variety as a factor.

\section{The Proof of the FL}

Shifting notation (as justified in \cite{Wald:2006},~\cite{CHL:2010}),
we let $F$ be the field of rational functions on a curve $X$ over a
finite field $k$.  One of the novelties of Ng\^o's work is to treat
the FL as identities over the global field $F$, rather than as local
identities at a given place of $X$.
By viewing each global section of $\OO_X(D)$ as a
rational function on $X$, each point $a\in{\AA}$ is identified with an
$F$-valued point $a\in \cc(F)$.  The preimage of $a$ under the
characteristic map $\chi$ is a union of conjugacy classes in $\gg(F)$,
and therefore corresponds to terms of the Arthur-Selberg trace formula
for the Lie algebra.  The starting point of Ng\^o's work is the
following geometric interpretation of the trace formula.

\begin{thm}[Ng\^o] There is an explicit
  test function $f_D$, depending on the divisor $D$,
  such that for every anisotropic element $a\in\AA^{an}$, the
  sum of the orbital integrals with characteristic polynomial $a$  
  in the trace formula for $f_D$ 
  equals
  the number of Higgs pairs in the Hitchin fibration over $a$,
  counted with multiplicity.
\end{thm}

The proof is based on Weil's description of vector bundles on a curve
in terms of the cosets of a compact open subgroup of $G(\ring{A})$. 
Orbital integrals have a similar coset description.

From this starting point, the past thirty years of research
on the trace formula can be translated into geometrical properties of
the Hitchin fibration.  In particular, Ng\^o formulates and then solves
the FL as a statement about counting points in Hitchin
fibrations.  

The identities of the FL are between the orbital integrals on two
different reductive groups $G$ and $H$.  A root system is associated
to each reductive group.  There is a duality of every root system that
interchanges its long and short roots.  The two reductive
groups of the FL are related
only indirectly: the root system dual to that of $H$ is a subset of
the root system dual to that of $G$ (Figure~\ref{fig:svdw}).  Informally, the set
of representations of a group is in duality with the group itself, so
by a double duality, when the dual root systems are directly related,
we might also expect their representation theories to be directly
related.  This expectation is supported by an overwhelming amount
of evidence.

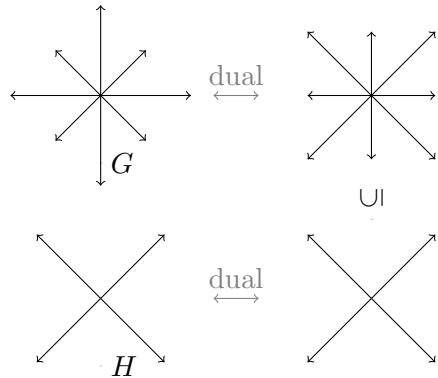
\tikzfig{svdw} {The two root systems in each row are in duality.  The
  root system on the bottom right is a subset of the root system on
  the upper right.}  
{
\begin{scope}[scale=0.03,xshift=-60cm]
\path (0:0) coordinate (P0);
\path (45:28.3) coordinate (P1) ;
\path (135:28.3)  coordinate (P2) ;
\path (225:28.3) coordinate (P3) ;
\path (315:28.3) coordinate (P4) ;
\path (0:40) coordinate (P5) ;
\path (90:40) coordinate (P6) ;
\path (180:40) coordinate (P7) ;
\path (270:40) coordinate (P8) ;
\foreach \i in {1,...,8}
{
  \draw[->] (P0) -- (P\i);
}
\fill (270:30) circle (0) node[anchor=west] {$G$};
\draw[<->,color=gray] (0:50) -- node[anchor=south] {dual} (0:70) ;
\end{scope}
\begin{scope}[scale=0.03,xshift=60cm]
\path (0:0) coordinate (P0);
\path (45:40) coordinate (P1) ;
\path (135:40)  coordinate (P2) ;
\path (225:40) coordinate (P3) ;
\path (315:40) coordinate (P4) ;
\path (0:28.3) coordinate (P5) ;
\path (90:28.3) coordinate (P6) ;
\path (180:28.3) coordinate (P7) ;
\path (270:28.3) coordinate (P8) ;
\foreach \i in {1,...,8}
{
  \draw[->] (P0) -- (P\i);
}
\fill(270:55) circle(0) node[anchor=west,rotate=90] {$\subseteq$}; 
\end{scope}
\begin{scope}[scale=0.03,xshift=-60cm,yshift=-90cm]
\path (0:0) coordinate (P0);
\path (45:40) coordinate (P1) ;
\path (135:40)  coordinate (P2) ;
\path (225:40) coordinate (P3) ;
\path (315:40) coordinate (P4) ;
\foreach \i in {1,...,4}
{
  \draw[->] (P0) -- (P\i);
}
\fill (270:30) circle (0) node[anchor=west] {$H$};
\draw[<->,color=gray] (0:50) --node[anchor=south] {dual} (0:70) ;
\end{scope}
\begin{scope}[scale=0.03,xshift=60cm,yshift=-90cm]
\path (0:0) coordinate (P0);
\path (45:40) coordinate (P1) ;
\path (135:40)  coordinate (P2) ;
\path (225:40) coordinate (P3) ;
\path (315:40) coordinate (P4) ;
\foreach \i in {1,...,4}
{
  \draw[->] (P0) -- (P\i);
}
\end{scope}
}

By using the same curve $X$ for both $H$ and $G$, and by comparing the
the characteristic maps for the two groups, Ng\^o produces a map
$\nu:\AA_H\to \AA_G$ of the bases of the two Hitchin fibrations, but
to kill unwanted monodromy he prefers to work with a base-change
$\tilde\nu:\tilde\AA_H\to \tilde\AA_G$.
The particular identities of the FL pick out a
subspace $\tilde\AA_\kappa$ of  $\tilde\AA_G$ 
containing $\tilde\nu(\tilde\AA_H)$.
Restricting the Hitchin
fibration to anisotropic elements, to prove
the FL, he must compare fibers of the two 
(base-changed, anisotropic) Hitchin
fibrations $\tilde\MM^{an}_G\to\tilde\AA^{an}_\kappa$ and
$\tilde\MM^{an}_H\to\tilde\AA^{an}_H$ over corresponding points of the base
spaces.

The base $\tilde\nu(\tilde\AA^{an}_H)$ contains a dense open subset of
elements that satisfy a transversality condition.  For
$\gg=\op{end}(V)$ this condition requires the self intersections of
the spectral curve (Equation \ref{eqn:spectral}) to be transversal
(Figure~\ref{fig:continuity}).  For a particularly nice open subset
$\tilde\UU\subset\tilde\nu(\tilde\AA^{an}_H)$ of transversal elements, the number of points in a
Hitchin fiber may be computed directly, and the FL can be verified in
this case without undue difficulty.
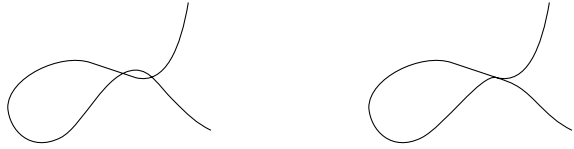
\tikzfig{continuity} {After giving a direct proof of the FL under the assumption
of transversality (left), Ng\^o obtains the
general case (right) by continuity.}
{
\begin{scope}[scale=0.20,xshift=-12cm]
\draw 
  (12.5,10) .. controls (12,7) and (11,4.5) .. (9,5)
   .. controls  (8,5.33) and (7,5.67) .. (6,6)
   .. controls (4,6.67) and (0.5,5) .. (0.5,3)
   .. controls (0.6, 1.5) and (2,0) .. (4,1)
   .. controls (5.5,1.75)   and (7,5.5) .. (9,5.5)
   .. controls (10,5.5) and (10.5,4.5 ) .. (11.5,3.5)
   .. controls (12.5,2.5) and (13,2) .. (14,1.5);
\end{scope}
\begin{scope}[scale=0.20,xshift=12cm]
\draw 
  (12.5,10) .. controls (12,7) and (11,4.5) .. (9,5)
   .. controls  (8,5.33) and (7,5.67) .. (6,6)
   .. controls (4,6.67) and (0.5,5) .. (0.5,3)
   .. controls (0.6, 1.5) and (2,0) .. (4,1)
   .. controls (5.5,1.75)   and (8,5.33) .. (9,5)
   .. controls (10,4.67) and (10.5,4.5 ) .. (11.5,3.5)
   .. controls (12.5,2.5) and (13,2) .. (14,1.5);
\end{scope}
}

To complete the proof, Ng\^o argues by continuity, that because the
identities of the FL hold on a dense open subset of
$\tilde\nu(\tilde\AA^{an}_H)$, the identities are also forced to hold on the
closure of the subset, even without transversality.  The justification
of this continuity principle is the deepest part of his work.

Through the legacy of Weil and Grothendieck, we know the number of
points on a variety (or even on a stack if you are brave enough) over
a finite field to be determined by the action of the Frobenius operator on
cohomology.  To cohomology we turn.  After translation into this
language, the FL takes the form of a desired equality of (the semisimplifications of) two perverse
sheaves over a common base space $\tilde\nu(\tilde\AA^{an}_H)$.  By the
BBDG decomposition theorem, over the algebraic closure of $k$, the
perverse sheaves break into  direct sums of simple terms, each given
as the intermediate extension of a local system on an open subset
$Z^0$ of its support $Z$~\cite{BBDG:1982}.  The decomposition theorem
already implies a weak continuity principle; each simple factor is
uniquely determined by its restriction to a dense open subset of its
support.  This weak continuity is not sufficient, because it does not
rule out the existence of supports $Z$ that are disjoint from the open
set of transverse elements.

To justify the continuity principle, Ng\^o shows that the support $Z$
of each of these sheaves lies in $\tilde\nu(\AA^a_H)$ and
intersects the open set $\tilde\UU$ of transverse
elements.  In rough terms, the continuity principle consists in
showing that every cohomology class can be pushed out into the open.
There are two parts to the argument: the cohomology class first is
pushed into the top degree cohomology, and then from there into the
open.  In the first part, the abelian varieties mentioned above enter
in a crucial way.  By taking cap product operations coming from the
abelian varieties, and using Poincar\'e duality, a nonzero cohomology
class produces nonzero class in the top degree cohomology of a Hitchin
fiber.  This part of his proof uses a stratification of the base of
the Hitchin fibration and a delicate inequality relating the dimension
of the abelian varieties to the codimension of the strata.

In the second part of the argument, a set of generators of the top
degree cohomology of the fiber is provided by the component group
$\pi_0$ of a Picard group that acts as symmetries on the fibers.
Recall that the two groups $G$ and $H$ are related only indirectly
through a duality of root systems.  At this step of the proof, a
duality is called for, and Ng\^o describes $\pi_0$ explicitly,
generalizing classical dualities of Kottwitz, Tate, and Nakayama in
class field theory.  With this dual description of the top cohomology,
he is able to transfer information about the support $Z$ on the
Hitchin fibration for $G$ to the Hitchin fibration on $H$ and deduce
the desired support and continuity theorems.  With continuity in hand,
the FL follows as described above.

Further accounts of Ng\^o's work and the proof of the FL appear in
~\cite{Nadler:2010}, \cite{Dat:2004}, \cite{Arthur:2010},
\cite{DN:2010},
~\cite{CHLaumon:2010},~\cite{Cass:2010},~\cite{NBC:report:2010}.

\section{Applications}

Only in the land of giants does the profound work of a Fields medalist get
called a lemma.  Its name reminds us nonetheless that the FL was never
intended as an end in itself.  A lemma it is.  Although proved only
recently, it has already been put to use as a step in the proofs of
the following major theorems in number theory.

\begin{enumerate}
\item Arthur's forthcoming classification of automorphic
  representations of classical groups~\cite{Arthur:2011}.
\item The calculation of the cohomology of Shimura
  varieties and their Galois representations~\cite{Morel:2010},~\cite{Shin:2010}.
\item The Sato-Tate Conjecture for elliptic curves over a totally real
  number field~\cite{BGHT:2010}.
\item Iwasawa's Main Conjecture for
  $GL(2)$~\cite{Skinner-Urban:2010},~\cite{Skinner:2010}.
\item The Birch and Swinnerton-Dyer Conjecture  for a
positive fraction of all elliptic curves over $\ring{Q}$~\cite{BS:2010}.
\end{enumerate}

The proof of the following recent theorem invokes the FL.  
It is striking that this result in pure
arithmetic ultimately relies on the Hitchin fibration, which
was originally introduced in the context of completely integrable
systems!

\begin{thm}[\cite{BGHT:2010}] 
  Let $n_p$ be the number of ways a prime $p$ can be expressed as a
  sum of twelve integers:
\[
n_p = \op{card}\,\{ (a_1,\ldots,a_{12})\in\ring{Z}^{12} \mid p = a_1^2 +\cdots +
  a_{12}^2 \}.
\]
Then the real number
\[
t_p = \frac{n_p -  8 (p^5 + 1)}{32 p^{5/2}}
\]
belongs to the interval $[-1,1]$, and as $p$ runs over all primes, the
numbers $t_p$ are distributed within that interval according to the
probability measure
\[
\frac{2}{\pi} \sqrt{1-t^2} \, dt.
\]
\end{thm}

\raggedright
\bibliographystyle{plain} 
\bibliography{/Users/thomashales/Desktop/googlecode/flyspeck/latex/bibliography/all}

\end{document}